# Modelling and Differential Costs Evaluation of a Two Echelons Repair System through Infinite Servers Nodes Queuing Networks

M. A. M. Ferreira[1]

Instituto Universitário de Lisboa (ISCTE – IUL)
ISTAR-IUL, BRU - IUL, Lisboa, Portugal

**Abstract**

After a short report of results on infinite servers' queues systems, focusing on its busy period, using networks of queues with infinite servers nodes a model is constructed to study a two echelons repair system. These repair systems may be useful, for instance, in the operation of a fleet of aircraft, of shipping or of trucks. Additionally, it is shown how this model may be used in the evaluation of differential costs.

## 1 Introduction

In a $M|G|\infty$ queue system:

-Customers arrive according to a Poisson process at rate $\lambda$,

-Each customer receives a service which length is a positive random variable with distribution function $G(.)$ and mean value $\alpha$,

-There are infinite servers, that is: when a customer arrives it always finds an available server,

-The service of a customer is independent the other customers' services and of the arrival process.

An important parameter is the traffic intensity, called $\rho$, being

$$\rho = \lambda\alpha \qquad (1.1).$$

The $M|G|\infty$ queue has neither losses nor waiting. Note that there is no queue in the ordinary sense of the word. For these systems it is not so important to study the population process as for other ones with losses or waiting, being much more interesting to study other processes as, for instance, the busy period. The busy period of a queue system begins when a customer arrives there, finding it empty, and ends when a customer leaves the system letting it empty. During the busy period there is always at least one customer in the system.

---

[1]manuel.ferreira@iscte.pt

A network of queues is a collection of nodes, arbitrarily connected by arcs, instantaneously traversed by costumers, and

-An arrival process is associated to each node,

-There is a commutation process that commands the different customers' paths.

Call $J$ the network number of nodes. If $J < \infty$, the nodes are numbered $1,2,\ldots,J$ and $U = \{1,2,\ldots,J\}$. The arrival processes may be the result of exogenous arrivals- from the outside of the collection- and of endogenous arrivals- from the other collection nodes. A network is open if any customer can enter it or leave it. Along this work open networks, with infinite servers in each node, and Poisson exogenous arrival rate $\lambda$ are considered. So

$$\Lambda = \begin{bmatrix} \lambda_1 \\ \lambda_2 \\ \vdots \\ \lambda_j \end{bmatrix} \quad (1.2)$$

is the network exogenous arrivals rate vector. The rate $\lambda_j$, is the exogenous arrival rate at node $j, j = 1,2,\ldots,J$ and

$$P = \begin{bmatrix} p_{11} & p_{12} & \cdots & p_{1J} \\ p_{21} & p_{22} & \cdots & p_{2J} \\ \vdots & \vdots & & \vdots \\ p_{J1} & p_{J2} & \cdots & p_{JJ} \end{bmatrix} \quad (1.3)$$

is the commutation process matrix, being $p_{jl}$ the probability of a customer, after ending its service at node $j$, go to node $l$, $j,l = 1,2,\ldots,J$. The probability $q_j = 1 - \sum_{l=1}^{J} p_{jl}$ is the probability that a customer leaves the network from node $j, j = 1,2,\ldots,J$. It is supposed that $P$ does not change with time and is independent of everything that is happening in the network. A network of this kind is equivalent to a $M|G|\infty$ system with Poisson process arrivals at rate $\lambda$, where each customer service time is its sojourn time in the network. This will be evidenced in this paper. But first some results about the $M|G|\infty$ systems busy period, that are also useful for these networks of queues evidently, are presented. Then a model is built, using these networks with infinite servers in each node, to study a two echelons repair system of a fleet of aircraft, of shipping or of trucks. The customers are the failures. And its service time is the time that goes from the instant at which they occur till the one at which they are completely repaired. The results referred above, about the $M|G|\infty$ queue busy period, allow the determination of some system performance measures. The theory is illustrated with a very simple and short numerical example and finally a differential costs analysis is presented. Some of this work is presented in Ferreira et al. [8,9].

## 2 Some Considerations on the M|G|∞ Queue Busy Period

Be $B$ the $M|G|\infty$ queue system busy period length. The mean value of $B$ is, see Takács [3], whatever is $G(.)$,

$$E[B] = \frac{e^\rho - 1}{\lambda} \qquad (2.1).$$

Calling $R(t)$ the mean number of busy periods that begin in $[0, t]$ ($t = 0$, the time origin, is the beginning of a busy period) see Ferreira[4],

$$e^{-\rho}(1 + \lambda t) \leq R(t) \leq 1 + \lambda t \qquad (2.2).$$

Let $N_B$ be the mean number of the customers served during a busy period in the $M|G|\infty$ queue systems. According to Ferreira [6],

-If $G(.)$ is exponential

$$N_B^M = e^\rho \qquad (2.3).$$

-For any other service distribution function

$$N_B \cong \frac{e^{\rho(\gamma_s^2+1)}(\rho(\gamma_s^2+1)+1)+\rho(\gamma_s^2+1)-1}{2\rho(\gamma_s^2+1)} \qquad (2.4)$$

where $\gamma_s$ is the $G(.)$ distribution function coefficient of variation.

## 3 The Sojourn Time of a Customer Laplace Transform in a Network of Queues with Infinite Servers in each Node and Poisson Exogenous Arrivals

Note that for these queue networks
   -The sojourn time of a costumer in each node is the service time, since there is not waiting,
   -The sojourn times of a customer in the various nodes are independent.

Using matrixes which form is suggested by (1.2) and (1.3) a simple formula is deduced to the sojourn time Laplace Transform, as a function of the service times in each node Laplace Transforms; see Ferreira[5] and Ferreira and Andrade [7].

Be $T$ the network sojourn time of a costumer and $S_j$ its service time at node $j$, $j = 1,2,\ldots J$. Be $G(t)$ and $G_j(t)$, $j = 1,2,\ldots,J$ the $T$ ad $S_j$, $j = 1,2,\ldots,J$, distribution functions respectively, being $\overline{G}(s)$ and $\overline{G}_j(s)$, $j = 1,2,\ldots,J$ the $T$ and $S_j$, $j = 1,2,\ldots,J$, Laplace Transforms, respectively.

If $\Lambda(s) = \begin{bmatrix} \lambda_1 \overline{G}_1(s) \\ \lambda_2 \overline{G}_2(s) \\ \vdots \\ \lambda_J \overline{G}_J(s) \end{bmatrix}$ and $P(s) = \begin{bmatrix} p_{11}\overline{G}_1(s) & p_{12}\overline{G}_2(s) & \cdots & p_{1J}\overline{G}_J(s) \\ p_{12}\overline{G}_1(s) & p_{22}\overline{G}_2(s) & \cdots & p_{2J}\overline{G}_J(s) \\ \vdots & \vdots & & \vdots \\ p_{J1}\overline{G}_1(s) & p_{J2}\overline{G}_2(2) & \cdots & p_{JJ}\overline{G}_J(s) \end{bmatrix}$

it results that
$$\overline{G}(s) = \sum_{n=0}^{\infty} \lambda^{-1} \Lambda^T(s) P^n(s)(I - P)A \qquad (3.1)$$
being $A$ a column with $J$ 1's. The final formula is obtained noting that (3.3) may be put in the form

$$\overline{G}(s) = \lambda^{-1}\Lambda^T(s)\bigl(I - P(s)\bigr)^{-1}(I - P)A \qquad (3.2)$$

supposing that $|I - P(s)| \neq 0$.

## 4 A Two Echelons Repair System

Suppose a fleet of aircraft, of shipping or of trucks which failures repairs occur in a base or in a remote station. The whole failures detected in the base are repaired there. Some of the failures detected in the station are repaired in the base with probability $p$, being necessary to transport them to the base, and the others in the station. Here the service time is the time that goes from the instant at which the failure occurs till the one at which it is completely repaired. When it is necessary to transport an item with a failure from the remote station to the base it is assumed that it is immediately possible, being the service time, now, the time that the transport lasts. It is also supposed that the failures occur according to a Poisson process at rate $\lambda$, being some detected in the remote station with probability $q$ and the others in the base. This situation may be modeled as a network of queues with three nodes, see Fig.1, at which 1 is the base, 2 is the remote station and 3 considers the required transports from the remote station to the base. Representing the variables related with each node for the same letter, as in the former sections, but with an index related to the node, obviously, $\Lambda = \begin{bmatrix} \lambda_1 \\ \lambda_2 \\ \lambda_J \end{bmatrix} = \begin{bmatrix} (1-q)\lambda \\ (1-p)q\lambda \\ pq\lambda \end{bmatrix}$ and $P = \begin{bmatrix} 0 & 0 & 0 \\ 0 & 0 & 0 \\ 1 & 0 & 0 \end{bmatrix}$.

Being the system, globally, a $M|G|\infty$ queue and, after (3.2),

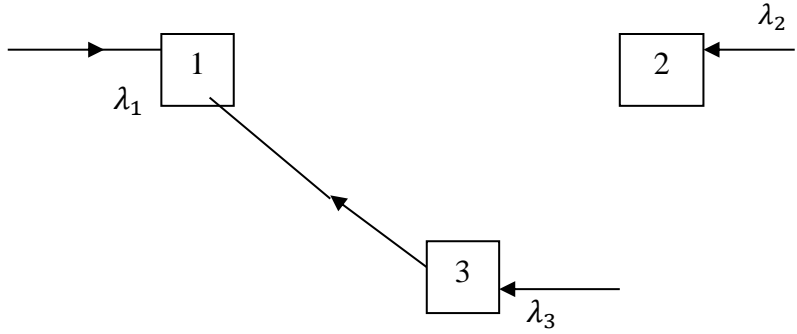

**Fig. 1.** The model network of queues in scheme

$\overline{G}(s) = [(1-q)\overline{G}_1(s) + (1-p)q\overline{G}_2(s) + pq\overline{G}_1(s)\overline{G}_3(s)]$. So, the service time distribution function is

$$G(t) = (1-q)G_1(t) + (1-p)qG_2(t) + pqG_{13}(t) \quad (4.1).$$

$G_{13}$ represents the distribution function of the convolution of the service time distributions in nodes 1 and 3. Three $M|G|\infty$ queues can also be considered. One related with the repairs in the base, which failures were detected there: $\lambda_b = (1-q)\lambda$, $G_b(t) = G_1(t)$ and $\rho_b = (1-q)\lambda\alpha_1$. Other related with the repairs in the remote station: $\lambda_{st} = (1-p)q\lambda$, $G_{st}(t) = G_2(t)$ and $\rho_{st} = (1-p)q\lambda\alpha_2$. And still other related with the repairs in the base after transport from the remote station: $\lambda_{tr} = pq\lambda$, $G_{tr}(t) = G_{13}(t)$ and $\rho_{tr} = pq\lambda(\alpha_1 + \alpha_3)$. Concerning the application of (2.4), being $\sigma_1^2$, $\sigma_2^2$ and $\sigma_3^2$ the variances corresponding to $G_1(.)$, $G_2(.)$ and $G_3(.)$, respectively $\gamma_{sb} = \gamma_{s1}$, $\gamma_{sst} = \gamma_{s2}$ and $\gamma_{str} = \frac{\sqrt{\sigma_1^2+\sigma_3^2}}{\alpha_1+\alpha_3}$. The coefficient of variation corresponding to the distribution function given for (4.1) is

$$\gamma_s = \sqrt{\frac{(1-q)(\sigma_1^2+\alpha_1^2)+(1-p)q(\sigma_2^2+\alpha_2^2)+pq(\sigma_1^2+\sigma_3^2+(\alpha_1+\alpha_3)^2)}{((1-q)\alpha_1+(1-p)q\alpha_2+pq(\alpha_1+\alpha_3))^2} - 1} \quad (4.2).$$

## 5 Application Example

Suppose that $G_1(.)$ and $G_2(.)$ are both exponential with mean 1 week and $G_3(.)$ is constant with value 1 week, $q = 0.3$ and $\lambda = \frac{1}{4weeks}$ (1 per month). Now $p = 0.9$ and considering 1 year (52 weeks) of operation is it possible to conclude if decreasing $p$ (that is: increasing the station capacity of repairs) there would be any advantage?

Making $p = 0.9; 0.8; \ldots 0.1$

-For the global system

Table 1. Global system

| $p$ | $\dfrac{e^\rho - 1}{\lambda}$ | $e^{-\rho}(1+\lambda 52)$ | $(1+\lambda 52)$ | $N_B$ |
|---|---|---|---|---|
| 0.9 | 1.51 | 10 | 14 | 2.04 |
| 0.8 | 1.45 | 10 | 14 | 2.03 |
| 0.7 | 1.40 | 10 | 14 | 2.02 |
| 0.6 | 1.40 | 10 | 14 | 2.03 |
| 0.5 | 1.35 | 10 | 14 | 2.02 |
| 0.4 | 1.29 | 11 | 14 | 2.01 |
| 0.3 | 1.24 | 11 | 14 | 1.99 |
| 0.2 | 1.24 | 11 | 14 | 2.00 |
| 0.1 | 1.19 | 11 | 14 | 1.99 |

-For the remote station

Table 2. Remote station

| $p$ | $\dfrac{e^{\rho_{st}} - 1}{\lambda_{st}}$ | $e^{-\rho_{st}}(1+\lambda_{st}52)$ | $(1+\lambda_{st}52)$ | $N_B$ |
|---|---|---|---|---|
| 0.9 | 1.00 | 1.38 | 1.39 | 1.01 |
| 0.8 | 1.01 | 1.75 | 1.78 | 1.02 |
| 0.7 | 1.01 | 2.12 | 2.17 | 1.02 |
| 0.6 | 1.02 | 2.48 | 2.56 | 1.03 |
| 0.5 | 1.02 | 2.84 | 2.95 | 1.04 |
| 0.4 | 1.02 | 3.19 | 3.34 | 1.05 |
| 0.3 | 1.03 | 3.54 | 3.73 | 1.05 |
| 0.2 | 1.03 | 3.88 | 4.12 | 1.06 |
| 0.1 | 1.03 | 4.22 | 4.51 | 1.07 |

- For the repairs in the base after transport

Table 3. Repairs in the base after transport

| $p$ | $\dfrac{e^{\rho_{tr}} - 1}{\lambda_{tr}}$ | $e^{\rho_{tr}}(1 + \lambda_{tr}52)$ | $(1 + \lambda_{tr}52)$ | $N_B$ |
|---|---|---|---|---|
| 0.9 | 2.14 | 3.94 | 4.51 | 3.49 |
| 0.8 | 2.12 | 3.65 | 4.12 | 3.13 |
| 0.7 | 2.11 | 3.36 | 3.73 | 2.81 |
| 0.6 | 2.09 | 3.05 | 3.34 | 2.54 |
| 0.5 | 2.08 | 2.74 | 2.95 | 2.30 |
| 0.4 | 2.06 | 2.41 | 2.56 | 2.10 |
| 0.3 | 2.05 | 2.07 | 2.17 | 1.91 |
| 0.2 | 2.03 | 1.73 | 1.78 | 1.76 |
| 0.1 | 2.02 | 1.37 | 1.39 | 1.62 |

- For the base

$\dfrac{e^{\rho_b}-1}{\lambda_b} = 1.09$ , $e^{-\rho_b}(1 + \lambda_b 52) = 8.48$ , $1 + \lambda_b 52 = 10.1$ and $N_B^M \cong 1.19$ , whatever is $p$.

Of course, in the operation of a fleet, it interests big idle periods and little busy periods. And if these occur it is good that they are as rare as possible, with a short number of failures. So, it is possible to conclude that the global system improves its performance as $p$ decreases but very lightly. The remote station grows worse has it was expected and the repairs in the base after transport improve its performance.

## 6 Differential Costs Evaluation

Note, for instance, that to guarantee less than two busy periods for the repairs in the base after transport, in the former example, it is essential that $p \leq 0.2$. But, for this, it is necessary to spend money in staff and material in the remote station. It is necessary to balance these expenses with the savings in transports. So, admitting that the transports costs are proportional to the probability $p$, call $c_i$ the cost associated to a certain initial value of $p$, $p_i$. If $p_i$ is decreased in a quantity $\Delta p_i$, the final cost $c_f$ associated to the new probability $p_i - \Delta p_i$ is related with $c_i$ as $c_f = c_i \dfrac{p_i - \Delta p_i}{p_i}$ being the differential cost $\Delta c_i = c_i - c_i \dfrac{p_i - \Delta p_i}{p_i} = c_i(1 - \dfrac{p_i - \Delta p_i}{p_i})$ that is

$$\Delta c_i = c_i \frac{\Delta p_i}{p_i} \quad (6.1).$$

The expression (6.1) gives the savings in transports. So, calling $k$ the investment in the remote station needed to decrease $p$, it is natural to request that $k \leq \Delta c_i$, that is

$$k \leq c_i \frac{\Delta p_i}{p_i} \quad (6.2)$$

to have no increase of financial problems. Note that all these costs must be considered as reported to a unit of time: 1 year, for example as it is usual.

## 7 Concluding Remarks

To apply this model, it is necessary to check if the failures occur according to a Poisson process. This hypothesis must be tested. After this it is very simple and correct to apply the other results. Carrillo[11] studied a model looking like the one presented here. But he did not consider either the possibility of transport from the station to the base or the busy period. So, these conceptions allow, in a very simple way, to evaluate the performance of this type of repair system as the example presented has shown. And permit also to evaluate, in a simple way, the financial viability of some type of modifications in the repair system.